\documentclass[a4paper]{article}
\usepackage{graphicx}
\usepackage{epsfig}
\usepackage{a4wide,epic}
\usepackage[cmex10]{amsmath}
\usepackage{amssymb}
\usepackage[caption=false,font=footnotesize]{subfig}
\usepackage{stfloats}
\usepackage{url}
\usepackage{ulem}
\usepackage[usenames,dvipsnames]{color}
\usepackage{pstricks}
\usepackage{enumitem}
\usepackage{rotating}
\setlist{nolistsep}

\newcommand{\un}[1]{\, \text{#1}}
\newcommand{\RR}{\mathbb{R}}
\newcommand{\CC}{\mathbb{C}}
\newcommand{\modal}[1]{\tilde{#1}}

\newcommand{\fns}[1]{\footnotesize{#1}}
\newcommand{\tightlist}{\setlength{\itemsep}{0cm} \setlength{\parskip}{0cm}}
\newcommand{\asedit}{}

\definecolor{gray}{rgb}{0.6,0.6,0.7}

\title{Control limitations from distributed sensing:\\ theory and Extremely Large Telescope application}
\author{A.~Sarlette and R.~Sepulchre\footnote{A.Sarlette is with SYSTeMS at Ghent University, Technologiepark Zwijnaarde 914, B-9052 Gent (Zwijnaarde), Belgium; \texttt{alain.sarlette@ugent.be}\newline\indent R.Sepulchre is with the Control Group at the University of Cambridge, Trumpington Street, Cambridge CB2 1PZ, UK.; and part-time with Systems and Modeling at University of Li\`ege, Department of Electrical Engineering and Computer Science (B28), B-4000 Li\`ege (Sart-Tilman), Belgium; \texttt{r.sepulchre@eng.cam.ac.uk}}}

\begin{document}
\normalem

\maketitle


\begin{abstract}
We investigate performance bounds for feedback control of distributed plants where the controller can be centralized (i.e.~it has access to measurements from the whole plant), but sensors only measure differences between neighboring subsystem outputs. Such ``distributed sensing'' can be a technological necessity in applications where system size exceeds accuracy requirements by many orders of magnitude. We formulate how distributed sensing generally limits feedback performance robust to measurement noise and to model uncertainty, without assuming any controller restrictions (among others, no ``distributed control'' restriction). A major practical consequence is the necessity to cut down integral action on some modes. We particularize the results to spatially invariant systems and finally illustrate implications of our developments for stabilizing the segmented primary mirror of the European Extremely Large Telescope.
\end{abstract}
\vspace{3mm}

\noindent {\footnotesize Keywords:\newline distributed detection, performance limitations, robustness, integral control, telescopes, distributed control.}


\section{Introduction}\label{sec:Intro}		

The massive availability of sensors and actuators in our environment is an opportunity to address large-scale problems by exploiting their interaction. The control of interacting localized subsystems (distributed plants) has drawn tremendous interest in the last decades, covering e.g.~agreement (consensus) in computer networks \cite{Tsitsiklis93}, synchronization of dynamical systems \cite{SujitSum,Synch} or collective robotic task solving \cite{BulloBook}. A defining property is the information sharing between subsystems --- information content, and interconnection topology. Common distinctions are reference-following \cite{BeardOnSats} vs.~autonomous coordination \cite{Synch} and centralized vs.~distributed control~\cite{B2Ptac11,Dahleh-Bamieh,Castro2002,Gorinevsky2008,Langbort2005,Stewart2003}. In centralized control, each local action is a function of measurements all over the system. 
The distributed control paradigm \cite{Dahleh-Bamieh,Castro2002,Langbort2005,Stewart2003,Gorinevsky2008} imposes a localized coupling in closed-loop: each local action depends on neighboring subsystem outputs only.

This paper shows how, ahead of the controller choice, structural restrictions on the \emph{sensing architecture} of a distributed plant can fundamentally constrain the performance of feedback. Specifically, we consider systems which only sense \emph{differences between neighboring subsystem outputs}. Unlike in distributed control, we allow the resulting measurements to be used in any --- in particular, centralized --- control computation. We call this ``local relative sensing'' or ``distributed sensing''. It is motivated by applications in which communication capabilities allow to quickly broadcast all measurements and control signals --- questioning a priori restrictions on controller structure --- but sensor technology does not allow accurate enough absolute measurements over the entire plant. This occurs in multi-scale problems, where accuracy requirements and plant size differ by many orders of magnitude.
The setting is inspired by our study of primary mirror stabilization for the European Extremely Large Telescope (EELT)~\cite{EELTreport}. We therefore propose indicative analytical results --- for a general case and for 1-degree-of-freedom spatially invariant systems --- followed by an illustration on this case study. We focus on two concerns. First, how distributed sensing influences the sensor noise vs.~disturbance rejection tradeoff, using the sensitivity transfer functions of classical linear control theory~\cite{Astrom2000,Astrom-Murray}. Second, how measurement model errors affect robustness. A major concrete consequence is the necessity to cut down integral action on some modes. 

Effects of noise and perturbations in distributed systems are examined in various ways in the literature. The authors of~\cite{B2Ptac11,Barooah1,Barooah2} restrict not only sensing, but also control, to local relative coupling (this is \emph{distributed control}). The authors of \cite{BarHesp1} study, in a static setting, bounds on the reconstruction of absolute position with respect to a leader from noisy local relative measurements. Optimal controllers for \emph{spatially invariant} plants are investigated in \cite{Dahleh-Bamieh}; for a locally coupled plant, the optimal gains decay exponentially as a function of distance between actuated and measured subsystems. This supports the use of distributed control, for which \cite{B2Ptac11} investigates performance limitations on a benchmark spatially invariant system. The robustness issue, that we first raised in \cite{EELTreport}, has been observed numerically with $\mu$-analysis for the segmented mirror application \cite{MacMyn2009}. Segmented mirror stabilization has been investigated by a few teams associated to Extremely Large Telescope projects \cite{SXJ2009,MacMartin2003,MacMyn2009}. In \cite{SXJ2009} the distributed sensing issue is put aside by assuming absolute measurements.

The paper is organized as follows. Section~\ref{sec:TheorA} formalizes distributed sensing and gives two motivating examples: a benchmark vehicle-chain problem and segmented mirror stabilization. Section~\ref{sec:TheorB} formulates how sensor noise (\ref{ssec:Noise}) and model errors (\ref{ssec:Robust}) induce performance limitations. Section~\ref{sec:TheorD} particularizes to 1-degree-of-freedom spatially invariant systems. Section~\ref{sec:EELT} illustrates our point on EELT primary mirror stabilization.


\noindent \textbf{Notation:} We write $i=\sqrt{-1}$ the imaginary unit. The element on row $j$, column $k$ of matrix $C \in \CC^{l \times m}$ is denoted $(C)_{j,k}$. $C^T$ and $C^*$ respectively denote transpose and complex conjugate transpose of $C$, and $\otimes$ the Kronecker product of two matrices. We denote $c \in \CC^l$ a column vector, $\Vert c \Vert = \sqrt{ \sum_k \vert c_k \vert^2 }$ its Euclidean norm. $I_m \in \RR^{m \times m}$ is the identity matrix and $\mathbf{1}_m \in \RR^m$ the column-vector of all ones. We interpret $s+C = s\, I_m + C$ if $s \in \mathbb{C}$ and $C \in \CC^{m \times m}$. For $D \in \CC^{m \times m}$ diagonal and $f$ a scalar function, $Y=f(D) \in \CC^{m \times m}$ is diagonal with $(Y)_{k,k} = f((D)_{k,k})$ for all $k$.


\section{Distributed sensing models}\label{sec:TheorA}


\begin{figure}[!h]
	\centering
	\setlength{\unitlength}{0.85mm}
	\begin{picture}(75,30)(-3,0)
    \put(0,8.5){\line(1,0){18.5}}
    \put(18,8){$\oplus$}
    \put(19,0){$n$}
    \put(10,10){$z$}
    \put(20,3){\vector(0,1){5}}
    \put(35,8.5){\vector(-1,0){14}}
    \put(35,5){\framebox(15,8){$B+\Delta$}}
    \put(30,1){\footnotesize \textsf{local relative map}}
    \put(75,8.5){\vector(-1,0){25}}
    \multiput(0,8.5)(75,0){2}{\line(0,1){12}}
    \put(0,20.5){\vector(1,0){7}}
    \put(7,17){\framebox(12,8){$C(s)$}}
    \put(19,20.5){\vector(1,0){19}}
    \put(37.8,19){$\oplus$}
    \put(39.7,26.3){\vector(0,-1){5}}
    \put(39,27){$d$}
    \put(40.5,20.5){\vector(1,0){5.5}}
    \put(46,17){\framebox(12,8){$G(s)$}}
    \put(47,27){\footnotesize \textsf{plant}}
    \put(5,27){\footnotesize \textsf{controller}}        
    \put(58,20.5){\line(1,0){17}}
    \put(27,22){$u$}
    \put(66,22){$y$}	
    \end{picture}
	\caption{Schematic representation of distributed sensing}\label{fig:modelA}
\end{figure}
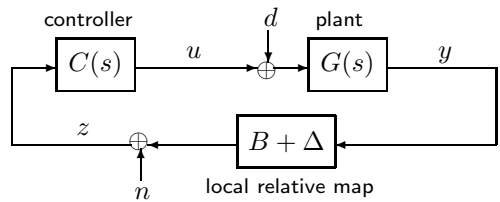

We consider a Laplace-domain model (see Fig.~\ref{fig:modelA})
\begin{eqnarray}
\label{eq:mod1}	y(s) & = & G(s) \; [\,u(s) + d(s)\,]\\
\label{eq:mod2} z(s) & = & [\, B + \Delta \,] \; y(s) + n(s)\\
\label{eq:mod3} u(s) & = & -C(s) \; z(s)
\end{eqnarray}
to represent $M \gg 1$ coupled $N$-dimensional subsystems. Components $kN$+$1$ to $(k$+$1)N$ of $y(s), u(s), d(s) \in \CC^{N_y}$ denote outputs, inputs and disturbances of subsystem $k$ in the Laplace domain, with $N_y$ = $NM$. We assume that the plant governed by $G(s)$ is stable. Output $z(s) \in \CC^{N_z}$ is obtained through the static map $[B+\Delta] \in \RR^{N_z \times N_y}$, where $B$ is the nominal sensor behavior and $\Delta$ a sensor model error. Each sensor measurement is corrupted by zero-mean independent identical Gaussian white noise, represented by $n(s)$ with covariance matrix $\sigma^2\, I_{N_z}$. For ease of presentation we assume $N_z \geq N_y$. The purpose of controller $C(s) \in \CC^{N_y\times N_z}$ is to reject disturbances $d(s)$ from $y$. Importantly, we do not restrict the controller (\ref{eq:mod3}) to be distributed, i.e.~we allow $C(s)$ to be a full matrix. We also allow the disturbances on different subsystems to be correlated, by investigating how a general vector $d(s)$ affects the controlled plant. This differs from e.g.~\cite{B2Ptac11,Barooah2} which examine $y$ for a given disturbance distribution (and controller).

The central element of our investigation is \emph{local relative measurement}. Let $q_k = \{kN$+$1, kN$+$2,...,(k$+$1)N\}$.\vspace{2mm}\newline
\textbf{Definition 1:} $\;B$ gives (unit-gain) \emph{relative measurements} between subsystem outputs if for each $l \in \{1,...,N_z\}$ there exist $q_j,q_k$ such that\vspace{-3mm}
\begin{eqnarray}
\nonumber 1. && (B)_{l,m} = 0 \;\;\; \text{for } m \notin q_j \cup q_k \\
\nonumber 2. && {\textstyle \sum_{m \in q_j} \;\;} \vert(B)_{l,m}\vert \;=\; {\textstyle \sum_{m \in q_j} \;\;} (B)_{l,m} \\
\nonumber && =\;  {\textstyle \sum_{m \in q_k} \;\;} \vert(B)_{l,m}\vert \;=\; -{\textstyle \sum_{m \in q_k} \;\;} (B)_{l,m} \; =\; 1 \; .
\end{eqnarray}
That means, each row $l$ of $B$ measures the difference between a convex combination of outputs of subsystem $j$ and a convex combination of outputs of subsystem $k$. For $N=1$, $B^T$ would be the \emph{oriented incidence matrix} of some graph $\Gamma_B$, where subsystems are nodes and sensors are edges; for $N>1$, $B^T$ has the interpretation of a generalized incidence matrix, with matrix-valued weights on each edge~\cite{BarHesp1}. The (generalized) Laplacian matrix of $\Gamma_B$ is $L=B^T B$. If $(L)_{l,m} \neq 0$ for some $l \in q_j$ and $m \in q_k$, then subsystems $j$ and $k$ are connected in $\Gamma_B$.
\vspace{2mm}\newline
\textbf{Definition 2:} A spatial structure $\mathcal{S}$ of dimension $\gamma \in \mathbb{N}$ associates a position $p(k) \in \mathbb{R}^{\gamma}$ to each subsystem $k$ such that $\Vert p(k)-p(l) \Vert \geq 1$ for $l \neq k$.\vspace{2mm}\newline
\textbf{Definition 3:} Given a spatial structure and a fixed spatial range $\rho \geq 1$, measurement map $B$ gives \emph{local relative measurements of range $\rho$} if it gives relative measurements and it only connects in $\Gamma_B$ subsystems for which $\Vert p(k)-p(l) \Vert \leq \rho$. We call this \emph{distributed sensing}.\newline
Many decentralized control settings associate a local measurement to each subsystem (graph node). With distributed sensing in contrast, \emph{measurements are the result of interactions between subsystems (graph edges).}

\textbf{Remark 1:} Local sensing has no meaning if it is not relative. Sensors giving ``absolute'' e.g.~positions actually physically measure positions with respect to a common (``central'') reference physically shared among all sensors. Absolute measurements thus correspond to centralized sensing. This is also acknowledged in the robotics community, distinguishing local$\cong$onboard from global$\cong$offboard sensors, see e.g.~\cite[Chapter 3]{BraunlBook}.

We study disturbance rejection limitations due to distributed sensing with $\tfrac{\rho}{M} \ll 1$. In the following two applications this arises as $M$ increases with the size of a large-scale plant while $\rho$ is limited by sensor technology.


\subsection{Vehicle chain application}\label{ssec:Cars}

A basic objective of e.g.~automated highway driving is to maintain constant inter-vehicle distance in a chain~\cite{SwaroopThesis}. Centralized sensing, typically the Global Positioning System (GPS), determines vehicle positions with respect to a common reference. GPS accuracy --- limited by atmospheric effects to a few meters --- is remarkable on the global scale, but likely insufficient to avoid collisions on a crowded highway. In local relative sensing, each vehicle directly measures the distance to its neighboring vehicles, easily up to centimeter-accuracy.

Define vehicle configuration $(y)_k = s_k - r_k \in \mathbb{R}$ with $s_k$ the coordinate of vehicle $k$ along the road and $r_k$ its desired coordinate, typically $r_k = k\, r_d$ in a moving frame with $r_d$ the desired distance between vehicle centers. Vehicles are controlled ($u$) and disturbed ($d$) by forces. Local relative sensors compare relative positions of consecutive vehicles, $(z)_k = s_{k+1}-s_k-r_d = (y)_{k+1} - (y)_k$, so $(B)_{k,k+1}=-(B)_{k,k}=1$ $\forall k$ and all other $(B)_{k,j}=0$ (path interconnection). This topology follows from a spatial structure $p(k) = k \in \mathbb{R}$ with $1<\rho<2$.

This benchmark problem has been studied before. String stability~\cite{SwaroopThesis} restricts its attention to a perturbation on the first vehicle and studies how it propagates along the chain, in absence of noise. \cite{B2Ptac11} imposes distributed \emph{control}, i.e.~each vehicle relies on local sensors only. Instead, we allow each vehicle to use information gathered by all sensors. By leaving the controller free, we also allow antisymmetric coupling, which is termed ``mistuned control'' in \cite{Barooah1} and improves disturbance rejection.


\subsection{Segmented mirror application}\label{ssec:Telescope}

\begin{figure}
    \centering
    \includegraphics[width=60mm]{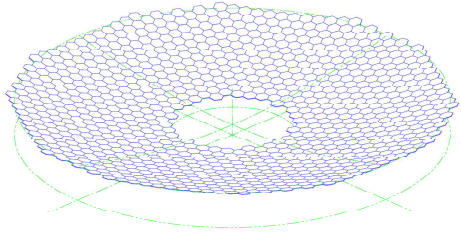}
    \caption{The EELT primary mirror made of 984 segments.}
    \label{fig:M1}
\end{figure}
\begin{figure} \centering
	\setlength{\unitlength}{0.85mm}
	\begin{picture}(90,30)(-3,0)
    \put(-5,0){\includegraphics[width=40mm]{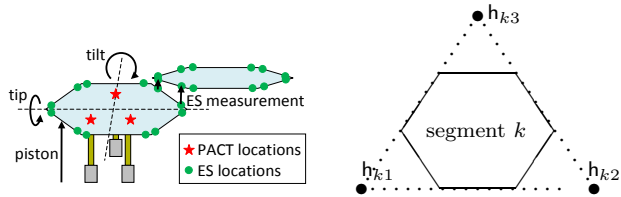}}
	\multiput(50,0)(1.8,0){18}{\circle*{0.5}}
	\multiput(50,0)(36,0){2}{\circle*{1.5}}
	\put(50,2){\fns{$\mathsf{h}_{k1}$}}
	\put(85,2){\fns{$\mathsf{h}_{k2}$}}
	\multiput(62.1,0.1)(0,17.9){2}{\line(1,0){11.8}}
	\put(68,27){\circle*{1.5}}
	\put(69.5,26.5){\fns{$\mathsf{h}_{k3}$}}
	\multiput(50,0)(1,1.5){18}{\circle*{0.5}}
	\multiput(86,0)(-1,1.5){18}{\circle*{0.5}}
	\put(62.1,0.1){\line(-2,3){5.93}}
	\put(73.9,0.1){\line(2,3){5.93}}
	\put(62.1,17.9){\line(-2,-3){5.93}}
	\put(73.9,17.9){\line(2,-3){5.93}}
	\put(60,8){\fns{segment $k$}}
	\end{picture}
    \caption{Left: actuation and sensing architecture. Position actuators (PACT) operate the segment's piston, tip and tilt degrees of freedom. Edge sensors (ES) measure local relative displacement of the segments. Right: points $(\mathsf{h}_{k1},\mathsf{h}_{k2},\mathsf{h}_{k3})$ used to define segment configuration.}
    \label{fig:segment}
\end{figure}

This paper has been inspired by our involvement \cite{EELTreport,EELTlaunch} in designing a controller to stabilize the segmented primary mirror of the European Extremely Large Telescope (EELT, see~\cite{Gilmozzi2007,ESOweb}). This telescope, run by the European Organization for Astronomical Research in the Southern Hemisphere (ESO), will offer unprecedented optical observation capabilities thanks to its primary mirror (M1) of world-record $42\un{m}$ diameter~\cite{Gilmozzi2007,ESOweb}. Its construction has started in 2012.
The scientific objective requires the light wave reflected by M1 to differ from an ideal one by less than $10$~nm root-mean-squared, after a linear filter that models adaptive optics corrections elsewhere in the telescope\footnote{See \cite{EELTreport,EELTlaunch} for more about this ``wavefront error''.}. One-piece mirrors meeting this accuracy under disturbances are currently limited to $\sim 9 \un{m}$ diameters. M1 is therefore composed of $M=984$ hexagonal segments of $0.7\un{m}$ edge length, see Fig.~\ref{fig:M1}, which are actively stabilized. The huge scale factor between the $10$~nm accuracy requirement and the $42$~m mirror size makes it technologically impossible to rely on segment displacement measurements with respect to a common reference. $N_z=5604$ ``edge sensors'' (ES) therefore measure the relative displacement, perpendicular to mirror surface, of adjacent segment edges (see Fig.~\ref{fig:segment}).

Disregarding in-plane motions of the segments (which have no optical effect), we consider $N=3$ degrees of freedom corresponding to piston, tip and tilt (PTT) of each segment i.e.~for a telescope pointing to zenith: vertical position and rotation around two horizontal axes; so $y \in \RR^{2952}$. We write the model after a coordinate change from the nominal mirror shape to a horizontal reference plane. Denoting $\pi_k$ the plane in $\RR^3$ which contains mirror segment $k$, define $(y_{3k-2},y_{3k-1},y_{3k})$ the height w.r.t.~reference plane of 3 points $\mathsf{h}_{k1},\mathsf{h}_{k2},\mathsf{h}_{k3} \in \pi_k$ positioned as depicted on Fig.~\ref{fig:segment}, right. Then the height of each point on segment $k$ is a convex combination of $(y_{3k-2},y_{3k-1},y_{3k})$. Each edge sensor measures the difference between two such convex combinations and thus fits Definition 1 of relative sensing. The spatial structure can be defined by $p(k)$= position of segment $k$'s center in the hexagonal lattice, and $0.7\sqrt{3}\text{m}< \rho < 0.7\cdot 3\text{m}$.

To control its 3 degrees of freedom, each segment is supported by 3 position actuators (PACT) which move perpendicularly to the mirror surface and whose force commands make up $u \in \RR^{2952}$. The small displacements allow for a linear model $G(s)$, discussed in Section \ref{sec:EELT}. Wind force is the strongest varying part in $d$, with characteristic frequencies below a few $\tfrac{1}{10}$Hz, and is spatially correlated among the segments. A second main component of $d$ are quasi-static disturbances: thermal effects and gravity on the moving mirror induce deformations of very low temporal frequency but high amplitude ($\simeq 1\un{mm}$). Noise $n$ is modeled as white Gaussian with $1 \un{nm}/\sqrt{\text{Hz}}$ power distribution for each sensor.


\subsection{Geometrical considerations on distributed sensing}\label{ssec:Interpretation}

Relative sensing reflects signal space invariance: sensors are insensitive to a common deviation of all system outputs since $B \, y = B \, ( y + \alpha \mathbf{1}_{N_y})$, for any $\alpha \in \mathbb{R}$. This {\it signal space} invariance \cite{Sarlette2009a}, involving $y$, should not be confused with the {\it spatial} invariance of \cite{B2Ptac11,Dahleh-Bamieh}, involving translations along spatial index $k$.
\vspace{2mm}
\newline\emph{Local} relative sensing can be viewed as measuring a spatial derivative: $\frac{y_k - y_l}{\Vert p(k)-p(l) \Vert}$ is the standard Euler discretization of the derivative of $y$ in direction $p(k)-p(l)$, evaluated at $\tfrac{p(k)+p(l)}{2}$. The approximation holds for $y$ varying on characteristic spatial scales much larger than $\Vert p(k)-p(l) \Vert$. For large $M$, an analogy with PDEs can therefore give insight for feedback design, see e.g.~\cite{Barooah1,Barooah2,Sarlette2009a}, although a rigorous link is tricky to establish~\cite{CurtainCrit}. The strong sensitivity to particular model errors $\Delta$ which we highlight in Section~\ref{ssec:Robust}, is analog to the drastic changes in PDE properties under small perturbations that change the dominant spatial derivatives.


\section{Distributed sensing limits performance}\label{sec:TheorB}



\subsection{Sensitivity and spectral graph properties}\label{ssec:Noise}

Like in classical control theory, the sensitivity and the complementary sensitivity are key transfer functions to capture the performance limitations of the feedback system. 
We take $\Delta=0$ but $n \neq 0$ in (\ref{eq:mod2}). Since $B\, B^T \, B$ has the same range and kernel as $B$ and there is no incentive to assign any controller gain to pure noise, we can write $C(s) = C_a(s) \, B B^T =:  C_b(s) \, B^T$ in (\ref{eq:mod3}) such that the closed-loop system (\ref{eq:mod1})-(\ref{eq:mod3}) becomes ($s$ dropped)
\begin{equation}\label{eq:Ranged}
y = [\,I_{N_y}+G\,C_b\,L\,]^{-1} \; [\,G\,d \, + \, G\,C_b\,B^T\,n \,]
\end{equation}
where $L=B^T B$. Write eigendecomposition $L = Q\,\Lambda\,Q^T$ with $Q$ orthogonal and eigenvalues $\lambda_k := (\Lambda)_{k,k}$ ordered as $\lambda_1 \geq \lambda_2 \geq ... \geq \lambda_{N_y} \geq 0$. Defining $\modal{y} = Q^T y$ and $K_b(s) = Q^T\, G(s) \, C_b(s) \, Q$ in (\ref{eq:Ranged}) we get 
$\modal{y} = [\,I_{N_y}+K_b\,\Lambda\,]^{-1} \; [\,Q^T\,G\, d + K_b\,Q^T\,B^T\,n \,]$.
Since we impose no restriction on the controller, we can view $K_b(s)$ as a freely tunable matrix transfer function\footnote{Issues like pole-zero cancellation and controller realizability, not specific to distributed sensing, are thereby ignored. They reduce to their SISO counterpart if $Q^T G(s) Q$ is diagonal in some ($s$-independent) basis, which occurs e.g.~when $G(s)$ is a scalar multiple of the identity matrix.}. Moreover, disregarding (temporarily) the singular modes of $\Lambda$, view $K(s)=K_b(s) \Lambda \in \mathbb{C}^{N_y \times N_y}$ as a tunable loop transfer function. Then (\ref{eq:Ranged}) becomes
\begin{equation}\label{eq:Modes}
	\modal{y}(s) = [\,I_{N_y}+K(s) \,]^{-1} \, \delta(s) + [\,I_{N_y}+K(s) \,]^{-1} K(s) \,\nu(s)
\end{equation}
where $\nu(s) := \Lambda^{-1}\,Q^T\,B^T\,n(s) \in \mathbb{C}^{N_y}$ is rescaled noise and $\delta(s) := Q^T\,G(s)\, d(s)$ is the open-loop deviation of $\modal{y}(s)$ induced by $d(s)$.
In the single-input / single-output case, expression (\ref{eq:Modes}) formulates the classical tradeoff between sensitivity $S = \frac{1}{1+K(s)}$ and complementary sensitivity $T = \frac{K(s)}{1+K(s)} = 1-S$. The specificity of distributed sensing is how $\nu$ relates to $n$, namely how the spectrum of the Laplacian $L$ scales the effect of measurement noise.\vspace{2mm}\newline
\textbf{Lemma 1:} \emph{If $\,n(s)$ is distributed as a zero-mean Gaussian of covariance $\,\sigma^2\, I_{N_z}$, then $\nu(s)$ is distributed as a zero-mean Gaussian of covariance $\sigma^2\, \Lambda^{-1}$.}\vspace{2mm} \newline 
\textbf{Proof:} Write the singular value decomposition $B = U \Sigma V^T$ where $U \in \RR^{N_z \times N_y}$ has orthonormal columns, $V \in \RR^{N_y\times N_y}$ is orthogonal and $\Sigma \in \RR^{N_y \times N_y}$ diagonal. $L=B^TB$ allows to identify $V=Q$ and $\Sigma = \Lambda^{1/2}$. Thus $\nu = \Lambda^{-1}\,Q^T Q \Sigma U^T\, n = \Lambda^{-1/2}\,U^T\,n$. As $U^T n$ is an orthonormal projection of $n$ on some subspace, it has a zero-mean Gaussian distribution with covariance matrix $\sigma^2\, I_{N_y}$. Gain $1/\sqrt{\lambda_k}$ on $(U^T n)_{k}$ multiplies its $\sigma^2$-variance by $1/\lambda_k$. \hfill $\square$

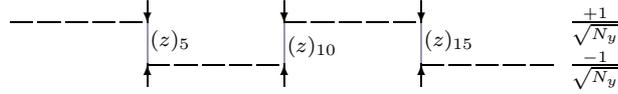
\begin{figure} \centering
	\setlength{\unitlength}{0.72mm}
	\begin{picture}(100,15)(5,0)
		\multiput(0,8)(5,0){5}{\line(1,0){4}}
		\multiput(25,0)(5,0){5}{\line(1,0){4}}
		\multiput(50,8)(5,0){5}{\line(1,0){4}}
		\multiput(75,0)(5,0){5}{\line(1,0){4}}
		\put(102,7){\footnotesize$\tfrac{+1}{\sqrt{N_y}}$}
		\put(102,-1){\footnotesize$\tfrac{-1}{\sqrt{N_y}}$}
		\multiput(25,0)(25,0){3}{\color{gray}\line(0,1){8}}
		\multiput(25,11.8)(25,0){3}{\vector(0,-1){4}}
		\multiput(25,-3.8)(25,0){3}{\vector(0,1){4}}		
		\put(25.5,3.5){\footnotesize$(z)_5$}
		\put(50.5,2.5){\footnotesize$(z)_{10}$}	
		\put(75.5,3.5){\footnotesize$(z)_{15}$}		
	\end{picture}
	\caption{Poorly observable configuration $y_i$ constructed for the proof of Proposition 1, in a 1D context with $N_y=20$, $N_c=4$, $\beta = 1/20$, assuming the sensing range extends to direct neighbors only. Each bar represents $(y_i)_k$ for a given $k$, with $y$ proportional to vertical position and subsystem index $k$ varying along the horizontal direction; $(z)_5$, $(z)_{10}$ and $(z)_{15}$ are the only measurements on the figure that differ from $0$.}\label{fig:NewFig}
\end{figure}

Unobservable modes, for which $\lambda_k = 0$, appear with infinite noise gain in (\ref{eq:Modes}) unless they have zero gain in $K(s)$ --- the only actual possibility. As $L\, \mathbf{1}_{N_y} = 0$, there is at least one unobservable mode. The following proposition shows how distributed sensing implies that, for $M \gg 1$, many other modes $k$ have $\lambda_k \ll 1$. 
\vspace{2mm} \newline
\textbf{Proposition 1:} \emph{Consider distributed sensing with a given spatial structure of dimension $\gamma$, range $\rho \ge 1$, and let $D_m$ the maximal number of sensors connected to a given subsystem. Then for  any (small) $c,\,N_c > 0$, there exists a (large enough) number of subsystems $M$ such that: $L=B^T B$ has at least $N_c$ eigenvalues $\lambda_k$ smaller than $c$, for any model of distributed sensing among $M$ subsystems according to $\rho$, $\gamma$, and $D_m$.}\vspace{2mm}\newline 
\textbf{Proof:} A few computations on eigenvalue bounds yield: if there are $N_c$ orthonormal vectors $y_i \in \RR^{N_y}$ for which $\Vert B \, y_i \Vert \leq b = \tfrac{c}{N_c}$, then $L=B^T B$ has at least $N_c$ eigenvalues smaller than $c$. We explicitly build such $y_i$.\newline
Given $N_c$, bounds on $D_m$, $\rho$ and $\gamma$ for distributed sensing, and any $\beta > 0$, the spatial structure ensures that for sufficiently large $N_y$, any compatible distributed sensing system can be partitioned into $N_c$ groups of subsystems such that:\vspace{-3mm}
\begin{itemize}
\item[(i)] each group contains $N_y/N_c$ outputs and
\item[(ii)] the maximum number of sensors connecting a subsystem of group $j$ with a subsystem of group $k \neq j$, over all group pairs $(j,k)$, is bounded by $\beta \, N_y$.
\end{itemize}\vspace{-3mm}
We build the $y_i$ by assigning the same value $\in \{\, \tfrac{+1}{\sqrt{N_y}},\, \tfrac{-1}{\sqrt{N_y}} \,\}$ to all the output components of a same group; see Fig.~\ref{fig:NewFig} for a schematic illustration. Without loss of generality, assume $N_c$ to be a power of $2$. Then $N_c$ orthogonal such $y_i$ can be built by associating $+$ or $-$ signs to the groups according to the elements of a Walsh code, as used in CDMA communication (see e.g.~\cite{CDMAex}). This code includes the vector $y_0$ where the output is a multiple of $\mathbf{1}_{N_y}$ (unobservable, $\Vert L \, y_0 \Vert = 0$). For the other $y_i$, it assigns $-$ signs to precisely $N_c/2$ groups. Then from (ii), $B\, y_i$ has at most $\beta N_y \, \tfrac{N_c(N_c-1)}{4}$ nonzero components, of value $\pm 2/\sqrt{N_y}$, such that $\Vert B \, y_i \Vert \leq \sqrt{\beta N_c (N_c-1)}$. Thus taking the initial $\beta \leq c^2 / (N_c^3 \, (N_c-1))$ yields the result. \hfill $\square$

The bound $D_m$ avoids that an arbitrarily large number of sensors can be used to perform the same measurement; this would indeed statistically improve the signal-to-noise ratio in a way that is not practically meaningful. The proof constructs poorly observable configurations by concerted deviations of large parts of the distributed system (e.g.~large subareas of the telescope mirror) inducing discontinuities that affect only a few sensors, see Fig.~\ref{fig:NewFig}. These are not the only ``poorly observable'' deformations. For settings with more structure, like in \cite{BarHesp1} or Section \ref{sec:TheorD}, one constructs even less observable deformations that have a different interpretation. Prop.1 focuses on worst-case modes, in contrast to \cite{B2Ptac11,Barooah2} which consider a mean over all modes. This is particularly relevant if, like in the EELT application, unfavorable $\delta$ are dominant. Indeed, for this application, a major disturbance source is wind force, which naturally features correlations on long spatial scales.



\subsection{Robustness to sensing model errors}\label{ssec:Robust}

We consider measurement model uncertainty $\Delta$ with
\begin{equation}\label{eqRob:Unc}
	(\Delta)_{k,l} \in [\,-\varepsilon\, (B)_{k,l} \,,\, +\varepsilon\, (B)_{k,l} \,] \quad \text{for all } k,l
\end{equation}
unknown, i.e.~each component of $B$ is (independently) subject to a relative uncertainty $\varepsilon$, for some given $\varepsilon \ll 1$. This uncertainty does not change the interconnection graph $\Gamma_B$.
Repeating the development of Section \ref{ssec:Noise} with $n=0$ and $\Delta \neq 0$ yields, similarly to \eqref{eq:Modes},
\begin{equation}\label{eqRob:Modes}
	\modal{y}(s) = [\,I_{N_y}+K(s) \,[I_{N_y}+\Phi] \,]^{-1} \, \delta(s)
\end{equation}
where $\Phi = \Lambda^{-1}\, Q^T \, B^T \,\Delta\, Q = \Lambda^{-1/2}\, U^T \Delta Q$, with singular value decomposition $B = U \Lambda^{1/2} Q^T$ as in Lemma 1. Note that in general $\Lambda$ is singular and the correct result is obtained by taking its pseudo-inverse in the definition of $\Phi$, i.e.~treating the lines $N_0+1$ to $N_y$ of $\Phi$ equal to zero, where $N_0$ is the index of the last nonzero eigenvalue of $\Lambda$. (By construction the columns $N_0+1$ to $N_y$ of $K$ equal zero.) For robust stability, the zeros of $\text{det}(\, I_{N_y}+K(s) \,[I_{N_y}+\Phi] \,)$ must have negative real part.

{\asedit We next analyze which properties on $\Phi$ are bad for stability. By standard matrix properties, $\text{det} (\, I_{N_y}+K(s) \,[I_{N_y}+\Phi] \,) = \text{det}(\, I_{N_0}+\bar{K}(s) \,[I_{N_0}+\bar{\Phi}] \,)$, where $\bar{K}$ (resp.$\bar{\Phi}$) contains the $N_0$ first columns and rows of $K$ (resp.$\Phi$). If further $K(s)$ is stable, we have for $s$ with positive real part: $\;\text{det}(I_{N_0}+\bar{K}(s)[I_{N_0}+\bar{\Phi}]) =0 \;$ $\Leftrightarrow$ $\; \;$  $\;\text{det}(\bar{K}(s)^{-1} + I_{N_0} +\bar{\Phi}) = 0 \;$. Consider e.g.~the hypothetical case of a disturbance $\Phi=\Phi_b$ defined by 
\begin{equation}\label{eq:ApproxPhi}
	(\Phi_b)_{k,l} = \; \left\{\begin{array}{ll} 0 & \text{for } (k,l) \neq (b,b) \\ \phi_b < 0 & \text{for } (k,l) = (b,b)  \end{array} \right. \; .
\end{equation}
We will approach this situation with an explicit construction later. Then we get:}
\vspace{2mm} \newline
\textbf{Proposition 2:} \emph{If there exists a real $s>0$ for which $\bar{K}(s)^{-1}$ exists and has entries bounded by $\varepsilon_b < 1/N_0$, and a sensing model error $\Delta$ can lead to a disturbance $\Phi = \Phi_b$ as in \eqref{eq:ApproxPhi} with $|\phi_b| > 1+N_0\,\varepsilon_b$, then the system is not robustly stable.}
\vspace{2mm}\newline 
\textbf{Proof:} {\asedit With $\Phi = \beta \Phi_b$ and $\beta$ varying from $0$ to $1$, the Gershgorin disks containing the eigenvalues of $\bar{K}(s)^{-1} + I_{N_0} +\bar{\Phi}$ all remain within the positive right-half plane, except for the row/column $b$ which moves from fully within positive to fully within negative right-half plane. By continuity, the corresponding eigenvalue must pass through zero, i.e.~$\text{det}(\, I_{N_y}+K(s) \,[I_{N_y}+\Phi_\beta] \,)$ has a zero with positive real $s$.} \hfill $\square$

More precise results can be given if reasonable structures are assumed for $K(s)$. For instance, if $K(s)$ is diagonal and positive for real positive $s$ (modal negative feedback) then robust stability requires $(\bar{K}(s))_{b,b} \;<\; 1\,/\,\vert 1 + \phi_b\vert$ for all real positive $s$.
{\asedit The point is that, given a nominal model with distributed sensing, a $\Phi_b$ with significant value of $\phi_b$ can be constructed from very small uncertainties in $\Delta$.} Indeed, the particular error $\forall k,l$\, :
\begin{equation}\label{eq:ZeMethod1}
 (\Delta)_{k,l} = -\text{sign}[(U)_{k,b}] \cdot \text{sign}[(Q)_{l,b}] \cdot \varepsilon\,\vert(B)_{k,l}\vert
\end{equation}
yields a maximally negative 
$(\Phi)_{b,b} = - \tfrac{\varepsilon}{\sqrt{\lambda_b}}\, \sum_{k,l}\; \vert (B)_{k,l}\vert \cdot \vert (U)_{k,b} \vert \cdot \vert (Q)_{l,b} \vert \, =:\phi_b$. {\asedit Here small values of $\lambda_b$ can lead to a strongly amplified model error, in particular a $\vert \phi_b \vert >1$ despite retaining $\varepsilon$-small uncertainty on each component of $\Delta$.} Numerical investigations on practical examples show that often with (adaptions of) \eqref{eq:ZeMethod1}, for some of the largest $b\leq N_0$ we can obtain: $(\Phi)_{b,b}$ largely dominates all the other elements of $\Phi$, so that $\Phi \approx \Phi_b$ given by \eqref{eq:ApproxPhi} for our purposes. {\asedit This shows that distributed sensing makes us dangerously close to the conditions of Proposition 2.\vspace{2mm}}
\newline {\asedit \textbf{Consequence:} \emph{From \eqref{eq:Modes}, a limitation on $K(s)$ as expressed by Prop.2 implies limited disturbance rejection. Most strikingly, integral control cannot be used for all modes since this would imply arbitrarily small $\bar{K}(s)^{-1}$ as $s$ approaches $0$. In modal control, poorly observable modes have poor static disturbance rejection as a tradeoff for robustness.}\vspace{2mm}}
\newline
Note that a disturbance like \eqref{eq:ZeMethod1} does not retain the symmetry owing to relative sensing, i.e.~$\Delta$
does not satisfy 
\begin{equation}\label{eqRob:gain}
	{\textstyle \sum_l} \; (\Delta)_{k,l} \, = \, 0 \quad \forall k \, .
\end{equation}
It is by breaking the symmetry of relative sensing that small errors can create large disturbances, analogously to changing the type of a PDE. In contrast, uncertainties that preserve the relative sensing symmetry (such as uncertainties in each sensor's gain) are easily shown to be much less detrimental to robust stability. However, many physical situations can plausibly lead to model uncertainties that violate \eqref{eqRob:gain}. In the EELT example (see Fig.~\ref{fig:segment}), each sensor consists of two parts fixed on adjacent segment edges and its measured value depends on the 3-dimensional relative motion of these parts. As a result of mirror curvature, a displacement $y \rightarrow y + \alpha \mathbf{1}_{N_y}$ with $\alpha>0$ brings the segments closer together in the mirror plane. This relative displacement along an unmodeled degree of freedom slightly affects the measurement, although nominally it should not.
Even with a flat mirror, mechanical constraints could induce systematic deformations of the segments or of their supporting structure, thus affecting the measurement through misalignments of sensor pairs even if the modeled $y$-difference does not change. Residual sensitivity to absolute output also seems plausible e.g.~for relative pressure or temperature sensors. All such errors can break the relative sensing symmetry and thereby severely limit the performance of the feedback system.


\section{Spatially invariant systems}\label{sec:TheorD}

Even though they do not require spatial invariance, the performance limitations described in the previous section take a simpler form under that additional hypothesis because they lead to insightful expressions in the spatiotemporal frequency domain. We refer the reader to  \cite{Dahleh-Bamieh,B2Ptac11,Gorinevsky2008,Castro2002} for system analysis of (LTSI) linear time-invariant and spatially invariant systems. 

The hypothesis of spatial invariance requires that all subsystems are equal and repeat an identical interconnection pattern. It \emph{restricts the controller}\footnote{Although, optimal control solutions for spatially invariant plants yield spatially invariant controllers \cite{Dahleh-Bamieh}.}, but it does not impose distributed control. We take $N = 1$ (1-degree-of-freedom subsystems) for simplicity, but we allow spatial invariance on a $\gamma$-dimensional toroidal structure. The system is then completely decoupled in spatiotemporal frequency domain, so we use the Fourier transform instead of the singular value decomposition; the real/orthogonal matrices of Section \ref{ssec:Noise} then become complex/unitary, with no added difficulty.

Our notation indexes subsystems in the toroidal spatial structure by $k=(k_1,k_2,...,k_\gamma) \in \mathcal{D} := \{1,...,M_1\} \times \{1,...,M_2\} \times ... \times \{1,...,M_\gamma\}$ with $M_1 \, M_2 \,...\, M_\gamma = M$. Localized measurement allows subsystems $j$ and $k$ to be connected through $B$ only if their indices are close, e.g.~if $\vert (j_a-k_a)\text{mod}(M_a) \vert \leq \rho'$ for all $a\in \{ 1,2,...,\gamma\}$, with $\rho'$ related to $\rho$. We denote $\; \mathcal{L} = \{\; l \in \mathcal{D} \; : \; \exists \text{ sensor comparing subsystems } k \text{ and } (k+l)\text{mod}(M_1,...,M_\gamma)\;,\; \forall k \; \}$. We write $D$ the number of sensors involving one given agent, such that e.g.~$N_z/M=D/2$.

Under spatial Fourier transform, the closed-loop equation (\ref{eq:Modes}) decouples into a set of (spatial) modes
\begin{equation}\label{eq:modcahat}
\modal{y}_\xi(s) = \frac{1}{1 + K_\xi(s)} \, \delta_\xi + \frac{K_\xi(s)}{1 + K_\xi(s)} \, \nu_\xi \, ,
\end{equation}
indexed by spatial frequency $\xi=(\xi_1,\xi_2,...,\xi_\gamma)$, with $\xi_j \in \{\, \tfrac{2k\pi}{M_j} : k=0,1,...,M_j-1 \, \}$. The following result holds, assuming that each sensor measures the difference between two output values; i.e.~each row of $B$ contains one element $+1$, one element $-1$ and the rest zeros. 
\vspace{2mm} \newline
\textbf{Lemma 2:} \emph{If $\,n(s)$ is distributed as a zero-mean Gaussian of covariance $\,\sigma^2\, I_{N_z}$, then $\nu_\xi(s)$ is distributed as a zero-mean Gaussian of variance $\; \sigma^2 / \lambda_\xi \;$ where
\begin{equation}\label{eq:SInoiseNu}
\; \lambda_\xi = 2\; {\textstyle \, \sum_{l \in \mathcal{L}}} \; \sin(l^T\, \xi/2)^2 \;\;  \text{ for each } \xi\,.
\end{equation}
For any local interconnection structure, given $c \ll 1$ there are at least $w_c := \Pi_{j=1}^\gamma \, (2\text{\emph{floor}}(\tfrac{c\,M_j}{\sqrt{2D}\gamma\rho'\pi})+1)$ frequencies $\xi$ for which $\nu_\xi$ has variance at least $\sigma^2/c^2$.}
\vspace{2mm}
\newline \textbf{Proof.} Lemma 1 remains valid with the Fourier transform replacing the SVD. Then we have $\lambda_\xi = B_\xi^* \, B_\xi$, where $B_\xi \in \CC^{D/2}$ is the $\xi$-component of the spatial Fourier transform of $B$. Computing the latter analytically yields \eqref{eq:SInoiseNu}. Bound $w_c$, which is not tight for $\gamma > 1$, is obtained by counting all $\xi=(\tfrac{2k_1\pi}{M_1},...,\tfrac{2k_\gamma\pi}{M_\gamma})$ for which $\tfrac{k_j}{M_j} \notin (c_1,\,1-c_1) \; \forall j$, then redefining $c=\sqrt{2D} \gamma \rho' \pi c_1$.
\hfill $\square$

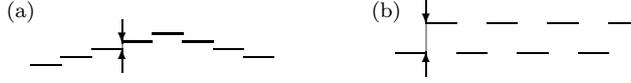
\begin{figure}[!t] \centering
\setlength{\unitlength}{0.8mm}
\begin{picture}(90,12)(2.5,0)
	\put(-4,8){\footnotesize (a)}
	\put(56,8){\footnotesize (b)}
	\multiput(0,0)(5,1.3){5}{\line(1,0){5}}
	\multiput(25,3.9)(5,-1.3){3}{\line(1,0){5}}
	\multiput(60,2)(10,0){4}{\line(1,0){5}}
	\multiput(65,7)(10,0){4}{\line(1,0){5}}
	\put(65,9){\color{gray}\line(0,-1){8}}
	\put(65,10.8){\vector(0,-1){4}}
	\put(65,-1.8){\vector(0,1){4}}
	\put(15,7){\color{gray}\line(0,-1){6}}
	\put(15,-1.2){\vector(0,1){4}}
	\put(15,7.6){\vector(0,-1){4}}
\end{picture}
\caption{Schematic representation of eigenvectors of $L$ for a 1-dimensional SI plant, similarly to Fig.~\ref{fig:NewFig}: (a) $\xi \approx 0$, long-range deformation; (b) $\xi=\pi$, short-range.}\label{fig:LongvsShortrange}
\end{figure}

The detrimental low $\lambda_\xi$ occur for frequencies $\xi$ close to $0$ modulo $2\pi$, i.e.~low spatial frequencies, corresponding to deformation modes of large characteristic length. Their poor observability can be understood as $B$ measuring a spatial derivative of the deformation \cite{Sarlette2009a}; see also Fig.~\ref{fig:LongvsShortrange}.

Regarding model errors, we can explicitly compute $\Phi$ if we restrict also $\Delta$ to be spatially invariant. This does not give as strong possibilities as for Proposition 2, but it still highlights a robustness issue.
\vspace{2mm} \newline
\textbf{Lemma 3:} \emph{For a spatially invariant uncertainty $\Delta$ as in \eqref{eqRob:Unc}, the possible $\Phi$ in an analog of \eqref{eqRob:Modes} diagonalized in spatial Fourier modes, are diagonal and include all $\beta \bar{\Phi}$ with $\beta \in [-1,1] \subset \RR$ and  
\begin{equation}\label{eq:Fourrob}
(\bar{\Phi})_{\xi,\xi} = \frac{\varepsilon \;\;{\textstyle \, \sum_{l \in \mathcal{L}}} \; \sin(l^T\, \xi/2)^2 \,\cdot\, \tfrac{i}{\tan(l^T\, \xi/2)}}{\phantom{\varepsilon \;}{\textstyle \, \sum_{l \in \mathcal{L}}} \; \sin(l^T\, \xi/2)^2 \phantom{ \,\cdot\, \tfrac{i}{\tan(l^T\, \xi/2)}}} \;\;  \text{ for each } \xi\,.
\end{equation}}\vspace{2mm} 
\newline \textbf{Proof.} We have $\Phi = \Lambda^{-1} [Q^* B^T (Q \otimes I)] \, [(Q^* \otimes I) \Delta Q]$ where $Q$ encodes spatial Fourier transform. Then $[Q^* B^T (Q \otimes I)]$ is a vertical concatenation of $D/2$ diagonal matrices: each one contains Fourier components $2i\, e^{-i l^T\xi/2}\, \sin(l^T \xi/2)$ of all $\xi$ for a given $l \in \mathcal{L}$ (avoid double counting). Taking `bad' $(\Delta)_{m,n} = \beta \varepsilon (B)_{m,n}$ for all $m,n$, we get $[(Q^* \otimes I) \Delta Q]$ a horizontal concatenation of $D/2$ diagonal matrices with Fourier components $2\, e^{i l^T\xi/2}\, \cos(l^T \xi/2)$.\hfill $\square$

\begin{figure} \centering
	\setlength{\unitlength}{1mm}
	\begin{picture}(95,65)
	\put(-8,35){\includegraphics[width=100mm]{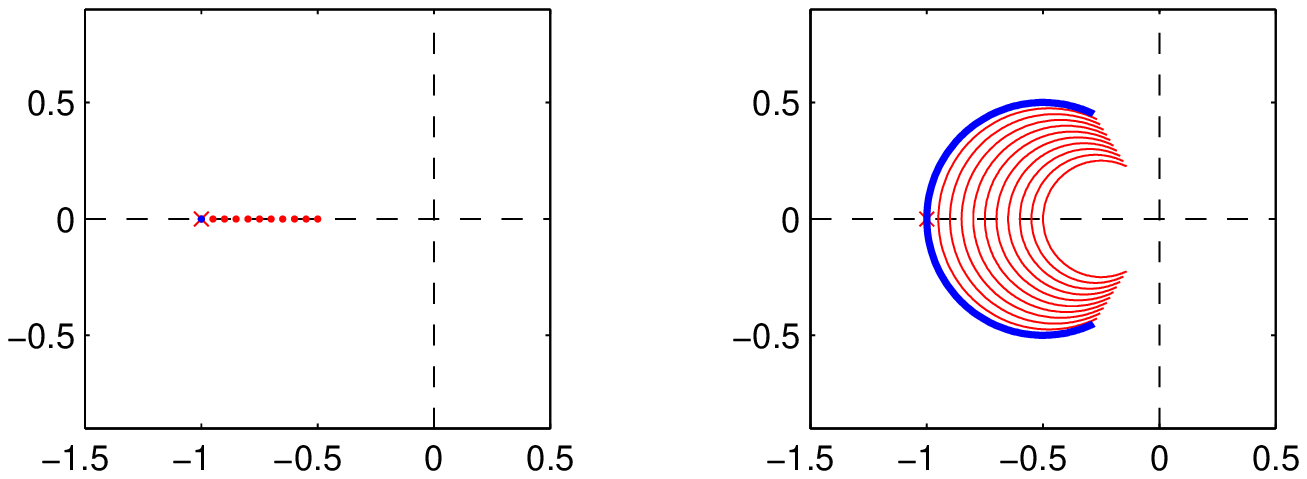}}	
	\put(-8,0){\includegraphics[width=100mm]{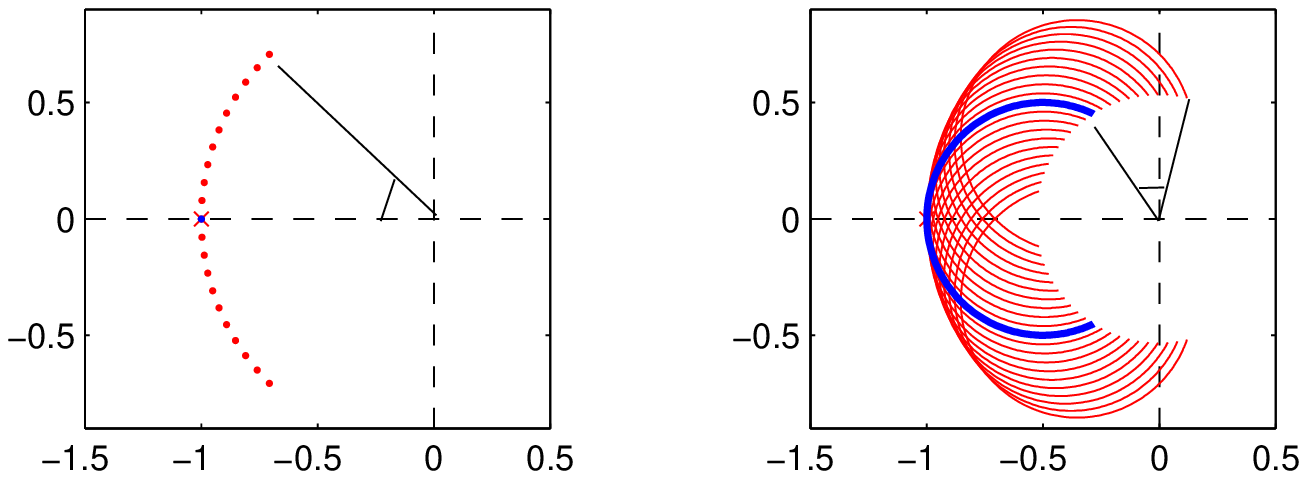}}
	\put(69.5,20){{\scriptsize$\phi_m$}}
	\put(21,17){{\scriptsize$\phi_m$}}
	\put(20,48){{\scriptsize$1/g_m$}}
	\put(22,52){\line(0,-1){2}}
	\end{picture}
	\caption{In dotted (left) and light (right) red: zones to exclude for Nyquist plot of $K_\xi(s)$ to have (top) gain margin $g_m$ and (bottom) phase margin $\phi_m$. Left: no/negligible sensing uncertainty, $K_\xi(s)$ must avoid circling $-1$. Right: ``bad'' sensing error $\beta \bar{\Phi}$ from Lemma 3, $K_\xi(s)$ must avoid circling an arc of circle (thick blue) that approaches $0 \in \CC$ up to $2\text{arctan}(1/\vert \bar{\Phi}_{\xi,\xi}\vert)$. [Plot values: $g_m=2$, $\phi_m=\pi/4$, $M=100$, $\varepsilon=0.05$ and $\xi=\tfrac{2\pi}{M}$, 1D nearest-neighbor coupling.]}\label{fig:PhaseMargin}
\end{figure}

Since all $l$ are restricted to small values by local sensing, for a low-frequency $\xi$ every term of the sum in the numerator is multiplied by a large factor $i\, / \, \tan(l^T\, \xi/2)$ with respect to the sum in the denominator, so the model error can have a dominant effect. Like for the general case, a `bad' $\Delta$ is one that changes the `meaning' of the measurement, i.e.~that violates \eqref{eqRob:gain}. `Bad' $\Delta$ have the following effect on robustness margins that account for other uncertainties in the system (e.g.~neglected dynamics, controller sampling / quantization, spatial invariance approximation).
\vspace{2mm} \newline
\textbf{Proposition 3:} \emph{Consider that robustness margins must be ensured under all spatially invariant measurement model errors $\Delta$ as in (\ref{eqRob:Unc}). Then the Nyquist plots of the nominal system ($\Delta=0$), at each spatial frequency $\xi$, must avoid enlarged ``exclusion zones'' as depicted on the right of Fig.\ref{fig:PhaseMargin}. In particular, for any given $\varepsilon,\kappa>0$ and fixed lattice dimension $\gamma$, there exist critical network sizes $M_{c,(i)}$ and $M_{c,(ii)}$ such that\newline
(i)  For all $M>M_{c,(i)}$, the set not to be encircled for nominal closed-loop stability, analog to $-1$ in a SISO system, has points in a $\kappa$-neighborhood of the origin for some $\xi \neq 0$.\newline
(ii) For all $M>M_{c,(ii)}$, the exclusion zone to ensure a phase margin $\phi_m=\kappa$ crosses the imaginary axis for some $\xi \neq 0$.}\vspace{2mm}
\newline \textbf{Proof.} Encirclement of $-1$ by the Nyquist plot of $K_\xi(1+\beta \bar{\Phi}_\xi)$ --- for $\bar{\Phi}$ as in Lemma 3 and $\beta \in [-1,1]$ --- is equivalent to encirclement by $K_\xi$ of the set $\mathcal{S}_\xi = \{ -1/(1+\beta (\bar{\Phi})_{\xi,\xi}) \; : \; \beta \in [-1,1] \}$. The latter set is an arc, extending by an angle $\theta=2\text{arctan}(\vert \bar{\Phi}_{\xi,\xi}\vert)$ on each side of the point $-1$, on a circle of radius $0.5$ centered at $-0.5 \in \CC$. The exclusion zone to robustly guarantee a gain (resp.~phase) margin $g_m$ (resp.~$\phi_m$) is obtained by rotating this arc around $0 \in \CC$ by all angles in $[-\phi_m,\phi_m]$ (resp.~scaling this arc towards $0$ by all ratios in $[1,g_m]$), see Fig.~\ref{fig:PhaseMargin}. A large enough system allows low enough $\xi$, such that $\vert \bar{\Phi}_{\xi,\xi}\vert$ defined by \eqref{eq:Fourrob} can get arbitrarily large, for any given $\varepsilon$. Then $\mathcal{S}_\xi$ tends towards a closed circle tangent to the imaginary axis.
\hfill $\square$

Proposition 3 shows how the spatially invariant measurement errors of Lemma 3, although not directly destabilizing, strongly diminish the robustness margins of an integral controller on low spatial frequencies $\xi$. This indicates that the system is particularly sensitive to dynamical model errors affecting those modes $\xi$: the combination of $\Delta$ with phase \& gain margins creates a kind of `trap' around the origin: the zone to avoid has an angular extension from $\sim \pi/2-\phi_m$ through $\pi$ to $\sim 3\pi/2+\phi_m$ and its inner radius can be arbitrarily close to $0$ for all angles in $[\tfrac{\pi}{2}\!-\!\phi_m,\,\tfrac{\pi}{2}\!+\!\phi_m] \cup [\tfrac{3\pi}{2}\!-\!\phi_m,\,\tfrac{3\pi}{2}\!+\!\phi_m]$. The following case study shows how, when spatial invariance is not imposed, worse situations as described in Section \ref{sec:TheorB} can indeed occur in practice.


\section{Case study: European Extremely Large Telescope primary mirror stabilization}\label{sec:EELT}

The ideas worked out in this paper are inspired by our study of the EELT primary mirror controller for and with ESO, which included review, proposal \& tuning of controllers, fault detection and isolation, specifications,... \cite{EELTreport}. We now show that the EELT system features the performance limitations studied above, and illustrate how they affect controller design. This supports the practical relevance of our observations. More details about the controller are given in the appendix.

\noindent\textbf{Model details:} We refer to Section \ref{ssec:Telescope} for a general modeling in the distributed sensing framework. The general model \eqref{eq:mod1} allows to include slight differences among segments and dynamical interactions through vibrations of the mirror-supporting structure (``backstructure''). Dynamics are however dominated by a static component plus mechanical oscillation of the actuators, such that we take $G(s)$ block-diagonal, with $3 \times 3$ blocks
\begin{equation}\label{eq:teldyns}
	G_0(s) =  ( \, s^2 J  + sb_a\, I_3 + k_a\, I_3 \, )^{-1} \, k
\end{equation}
where $k_a \in \RR$ is actuator stiffness, $b_a \in \RR$ gives (very small) dissipation, $J \in \RR^{3 \times 3}$ represents segment+actuator inertia matrix. Characteristic values for this mass-spring-damper subsystem are oscillation frequencies around $50$~Hz, with $\sim 1\%$ modal damping. The dominance of low frequency disturbances and the difficulty of high-frequency plant characterization (e.g.~backstructure vibrations) motivate a low-bandwidth model \& controller. We recall that the latter must seek an optimal tradeoff in presence of Gaussian white sensor noise $n$ with covariance matrix $\sigma^2 \, I_{5604}$ and $\sigma = 1$~nm/$\sqrt{\text{Hz}}$. An ES model uncertainty of the form \eqref{eqRob:Unc} is expected, with $\varepsilon$ up to $0.01$ and possibly all-independent components (see paragraph after \eqref{eqRob:gain}). We recall the main disturbances in $d$:
\newline $\bullet$ Mounting/wear errors and deformations due to gravity acting differently as the mirror moves, are static but large (micro- to millimeters), and mostly expected to be spatially uncorrelated.
\newline $\bullet$ Wind force is concentrated on low temporal frequencies, but it features strong long-range correlations in space, i.e. predominantly low spatial frequency.

\begin{figure}[!t] \centering
    \centering
\setlength{\unitlength}{1mm}
\begin{picture}(55,33)
	\put(0,3){\includegraphics[bb=58 206 585 610,clip=true,width=42mm]{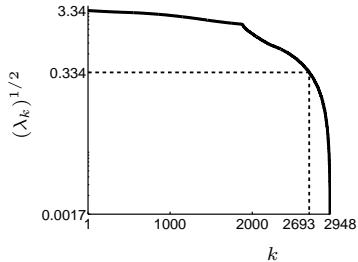}}
	\put(30,-1){\scriptsize $k$}
	\put(-4,15){\scriptsize \begin{sideways}$(\lambda_k)^{1/2}$\end{sideways}}
\end{picture}
    \caption{Singular values of $B$, i.e.~the diagonal elements $\lambda_k^{1/2}$ of $\Lambda^{1/2}$, in logarithmic scale. The last 4 modes are unobservable ($\lambda_k = 0$ for $k=2949$ to $2952$, beyond plot limits).}
    \label{fig:M1_singularvalues}
\end{figure}

\noindent\textbf{Limitations from distributed sensing:} Figure~\ref{fig:M1_singularvalues} shows the singular values $\lambda_k^{1/2}$ of $B$, in logarithmic scale. The last four modes are unobservable, $\lambda_k = 0$ for $k = 2949,...,2952$ (thus $N_0=2948$); they consist of full mirror translation, rotations and ``defocus'', and have to be addressed by other controllers based on wavefront sensing (see e.g.~\cite{EELTlaunch}). Among the observable modes, some indeed have very small $\lambda_k$ \textbf{[Section~\ref{ssec:Noise}, Prop.~1]} --- e.g.~the last 260 modes have a $\lambda_k < 0.01 \, \lambda_1$. This makes the performance limitations as expressed in Section \ref{ssec:Noise} relevant. The corresponding deformation modes (rows of $Q$, as $B^T B = Q \Lambda Q^T$) feature oscillations of long characteristic length, unfortunately where wind is strongest; the modes of the first $\lambda_k$ in contrast wildly oscillate from one segment to its neighbor. Regarding robustness, $\varepsilon = 0.01$ is insufficient to exclude all dangerous $\Delta$ \textbf{[Section~\ref{ssec:Robust}, Prop.~2]}: applying the procedure described in \eqref{eq:ZeMethod1} to $b=2948$, we get a $\Phi$ with one (exact) eigenvalue at $-9.7$. Thus large $K(s)$ --- e.g.~controllers with an integral term --- would bear the danger of instability, unless we can further reduce model uncertainty by a factor of at least $9.7$. This analysis agrees with the observations of~\cite{MacMyn2009}.

\noindent\textbf{General controller design:} We start from integral control $K_I/s$ to reject the expected low-frequency disturbances, but have to add integrator leakage $A_I$, yielding $(sI+A_I)^{-1} K_I$, to limit $\vert K(s) \vert$ for robustness. A double-pole is added to further damp the high-frequency resonance peaks due to the low damping of the actuator oscillator in \eqref{eq:teldyns}. This finally gives
\begin{equation}\label{eq:GeneralController}
	C(s) = -(sI+A_I)^{-1} K_I \; \frac{1}{(s/p+1)^2} \, .
\end{equation}
The last factor is a scalar, with fixed $p = 2 \pi\,20\un{rad/s}$. For $K_I,A_I$ we have analyzed two options: a centralized controller, based on diagonalization of the system into spatial eigenmodes \textbf{[Section \ref{sec:TheorB}]}; and a distributed controller, tuned in spatial frequency domain \textbf{[Section \ref{sec:TheorD}]}. The tradeoff between robust stability \textbf{[Section~\ref{ssec:Robust}]} and disturbance rejection \textbf{[Section~\ref{ssec:Noise}]} makes integrator leakage an essential tuning parameter.

\begin{figure}[!h]
    \centering
    \includegraphics[bb=20 180 585 610,clip=true,width=46mm]{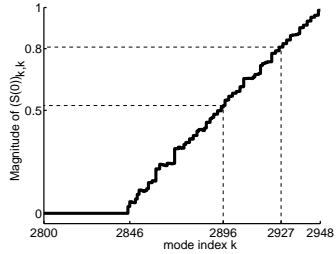}
    \caption{DC gain of disturbance-to-output transfer function, $\vert (S(0))_{k,k} \vert$, with centralized control in modal basis.}
    \label{fig:modal_SS_error}
\end{figure}

The centralized controller is analyzed in the \emph{modal basis}, obtained from the SVDecomposition of $B$ \textbf{[Section \ref{sec:TheorB}]} by $\modal{y} = Q^T y$ (``deformation modes''), $Q^T u,\, U^T z$. Both $K_I$ and $A_I$ are taken diagonal in this basis, writing their components $K_I(k)$ and $A_I(k)$, with $k=1,2,...,2952$ in order of decreasing $\lambda_k$. Due to the distributed sensing, small $\lambda_k$ --- which indeed correspond to mirror deformations of large characteristic length \textbf{[Lemma 2]} --- are sensitive both to noise \textbf{[Section~\ref{ssec:Noise}]} and to model uncertainties \textbf{[Section~\ref{ssec:Robust}]}, and this limits the control performance, especially the properties of its disturbance rejection sensitivity $S(s)=[I_{N_y}+K(s)]^{-1}$ at steady-state $s=0$ \textbf{[Eq.~\eqref{eq:Modes}]}. Figure~\ref{fig:modal_SS_error} shows the steady-state disturbance rejection $\vert (S(0))_{k,k} \vert \approx \left\vert \dfrac{1}{1+\sqrt{\lambda_k} K_I(k)/A_I(k)} \right\vert$ resulting from our robust modal tuning. For modes $k$ where leakage was not necessary $\vert (S(0))_{k,k} \vert = 0$; but for the modes of lowest $\lambda_k$ rejection of static errors becomes illusory \textbf{[Consequence of Prop.2 together with \eqref{eq:Modes}]}: e.g. $\vert (S(0))_{2948,2948} \vert \approx 1/(1+0.0017*5.7/0.7) = 0.986$. 

Because the performance limitations are independent from the controller design, one might wonder whether a distributed controller --- relying on local information only --- can attain the same performance as a centralized one. Benefits of distributed control include local communication requirements, arguably better robustness to failures, and a design method that builds directly on the physical distributed spatial structure, rather than on mathematically diagonalizing a large system matrix. We therefore design a controller where $(u)_k$ depends only on measurements made at the edges of segment $k$ and of its immediate neighbors. We further approximate the mirror as invariant  w.r.t.~translations from one segment to another (after all our basic model is, except at the mirror boundaries). This allows us to use the LTSI (Linear Time- and Space-Invariance) method described in \cite{Dahleh-Bamieh,Stein2005,Gorinevsky2008} for controller tuning in spatiotemporal frequency domain \textbf{[Section \ref{sec:TheorD}} extended to block-MIMO: one $3\times 3$ transfer matrix per spatial frequency].
\begin{figure} \centering
\setlength{\unitlength}{0.8mm}
\begin{picture}(90,12)(4,0)
	\put(-4,8){\footnotesize (a)}
	\multiput(0,0)(5,1.3){8}{\line(1,0){5}}
	\put(15,7){\color{gray}\line(0,-1){6}}
	\put(15,-1.2){\vector(0,1){4}}
	\put(15,7.6){\vector(0,-1){4}}
	\put(56,8){\footnotesize (b)}
	\put(75,7){\color{gray}\line(0,-1){6}}
	\put(75,-1.2){\vector(0,1){4}}
	\put(75,7.6){\vector(0,-1){4}}
	\multiput(60,3.9)(5,0){8}{\line(4,-1){5}}
\end{picture}
\caption{Spatial Fourier modes in a linear spatial structure with two degrees of freedom: piston and tilt. At low spatial frequency, a high-amplitude piston mode ($\mathbf{a}$) gives the same ES signal as a low-amplitude tilt mode ($\mathbf{b}$).}\label{fig:locestmodes}
\end{figure}
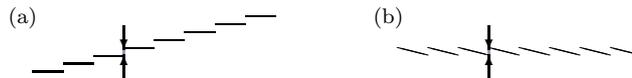
The performance limitations owing to distributed sensing are illustrated by considering how a local controller would estimate the overall mirror deformation from local measurements  \textbf{[Lemma 2]}. Indeed, the similar local sensor values can indicate global deformations of very different magnitude, depending on how they have to be extrapolated. Figure \ref{fig:locestmodes} depicts the most extreme case, where in a linear spatial structure, situations $\mathbf{a}$ and $\mathbf{b}$ would give the same ES signals, while the individual segments in $\mathbf{a}$ are displaced much more than in $\mathbf{b}$.

\begin{figure}[!h] \centering
	\setlength{\unitlength}{0.875mm}
	\begin{picture}(90,35)(0,35)
		\put(2,35){\includegraphics[width=37mm]{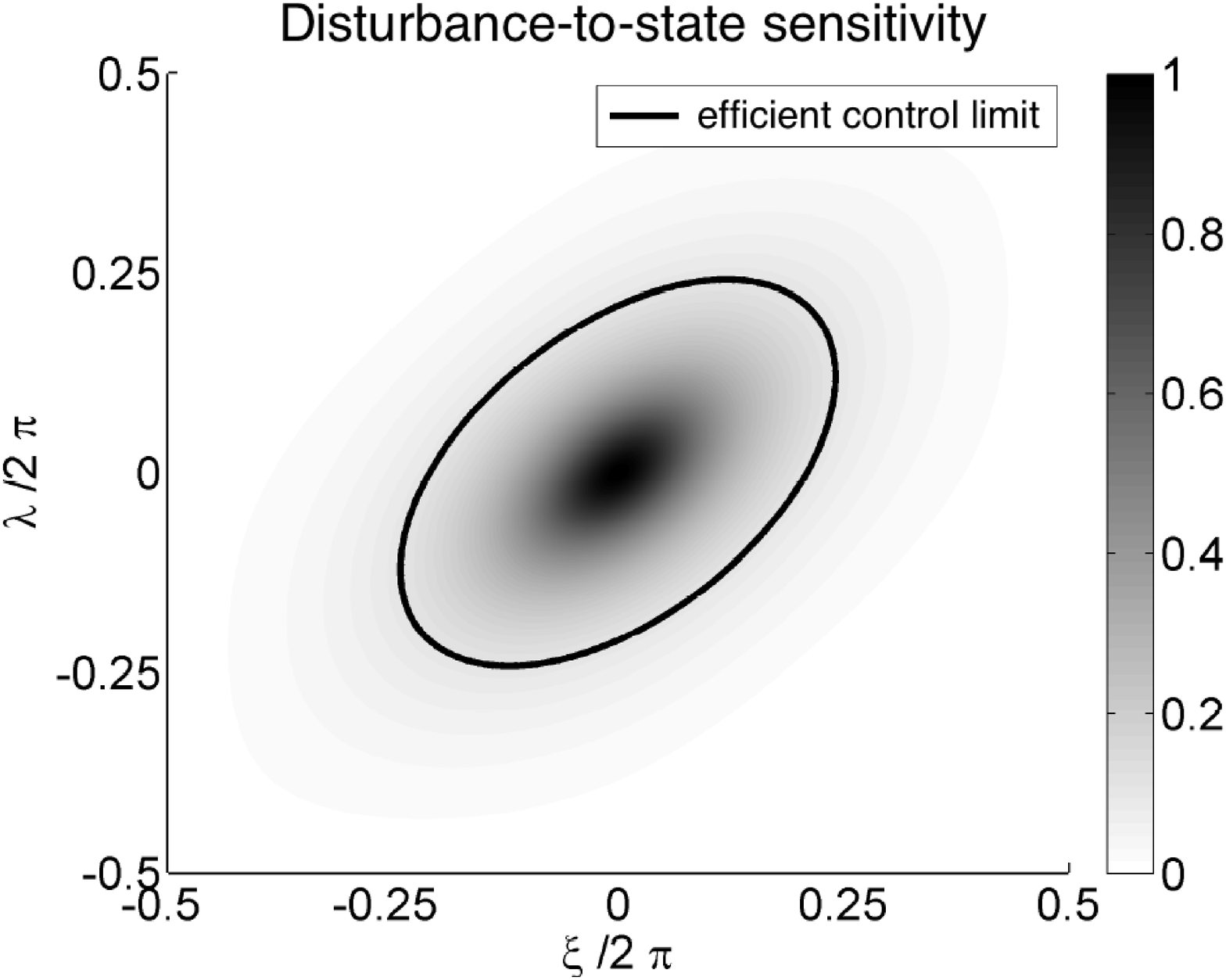}}
		\put(0,65){\footnotesize{a}}		
		\put(50,35){\includegraphics[width=37mm]{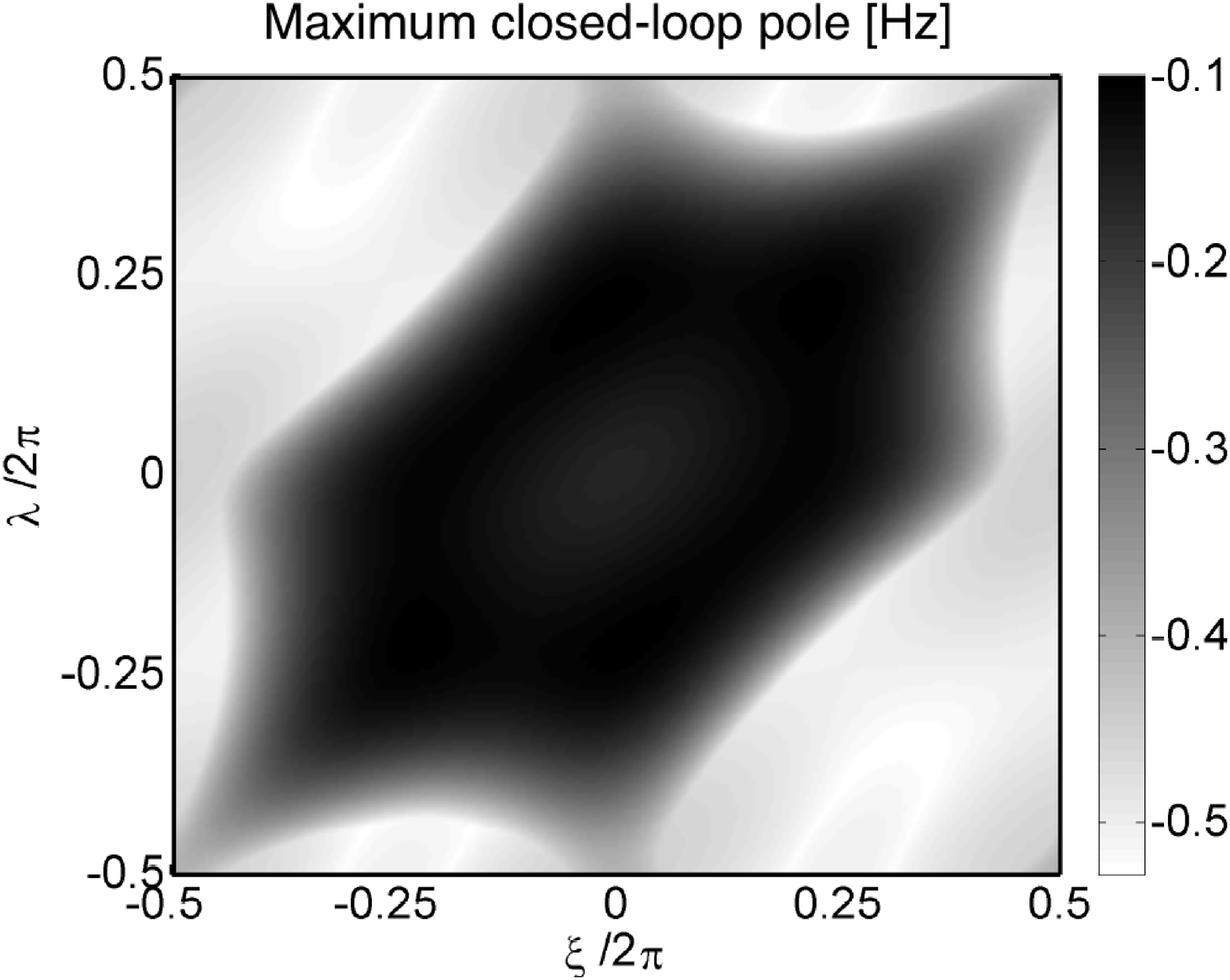}}
		\put(49,65){\footnotesize{b}}
	\end{picture}
    \caption{LTSI approximation properties of the distributed controller, on piston motion. See text for details.}
    \label{fig:LTSI_tuning_analysis}
\end{figure}

Properties of the distributed controller \emph{for an LTSI model} of the mirror \textbf{[Section \ref{sec:TheorD}]} are further illustrated on the steady-state disturbance rejection sensitivity $S(\xi;s\!\!=\!\!0) = [I_3+K_\xi(s\!\!=\!\!0)]^{-1}$ \textbf{[Eq.\eqref{eq:modcahat}} extended to block-MIMO], the LTSI equivalent of Fig.\ref{fig:modal_SS_error}. Tip and tilt Fourier modes are all sufficiently observable in the two-dimensional structure to allow zero leakage; lower observability at high spatial frequencies however requires to limit $K_I^*$. For piston, $K_I^*$ is limited at low to moderate spatial frequencies and, most notably, leakage is necessary at low spatial frequencies \textbf{[Lemma 3]} and induces large $S$, as shown on Fig.\ref{fig:LTSI_tuning_analysis}a. Figure~\ref{fig:LTSI_tuning_analysis}b shows the maximum real part of closed-loop poles of the MIMO system with uncertainties. The tradeoff resulting from distributed sensing is visible: at low spatial frequencies, where leakage limits $S$-performance, the closed-loop poles are pushed against the robustness limit. 

\noindent\textbf{Simulations:} We have simulated the controllers on a model of the full telescope composed of 20000 states extracted from a comprehensive finite element model. The model is subjected to a representative wind profile and to $1 \un{nm}/\sqrt{\text{Hz}}$ sensor noise. The limitations induced by distributed sensing \textbf{[Section \ref{sec:TheorB}]} do show up. When setting leakage $A_I$ to zero, simulations with `bad' $\Delta$ constructed from \eqref{eq:ZeMethod1} indeed feature unstable behavior, both for centralized \textbf{[Prop.2]} and distributed controllers \textbf{[Lemma 3]}. {\asedit And on the nominal system (with $\Delta=0$), choosing a constant closed-loop gain $K_I(i) \sqrt{\lambda_i}$ over all modes causes so much noise amplification that closed-loop performance is worse than in open-loop. This illustrates the reality of distributed sensing performance limitations.}

\begin{figure*}[!t]
\setlength{\unitlength}{0.8mm}
\begin{picture}(180,45)(-3,0)
	\put(-3,-3){\includegraphics[height=40mm]{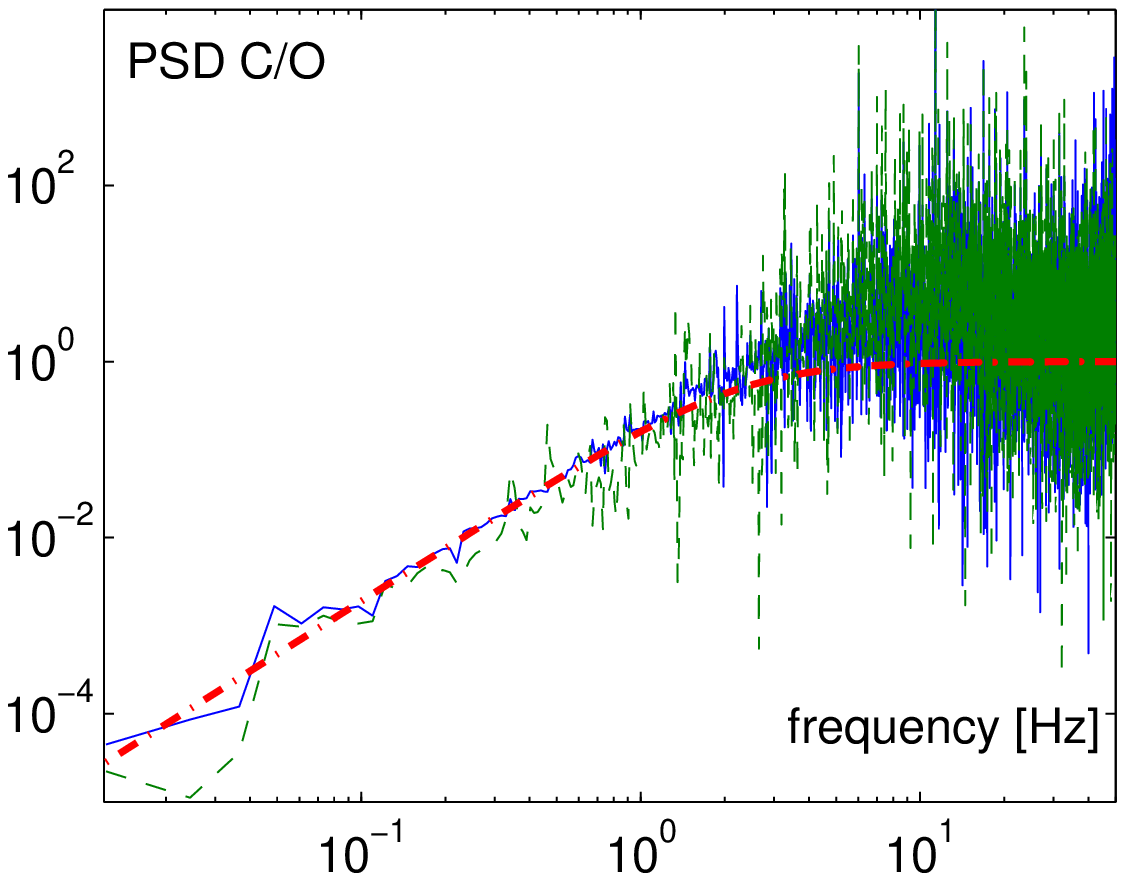}}
	\put(62,-3){\includegraphics[height=40mm]{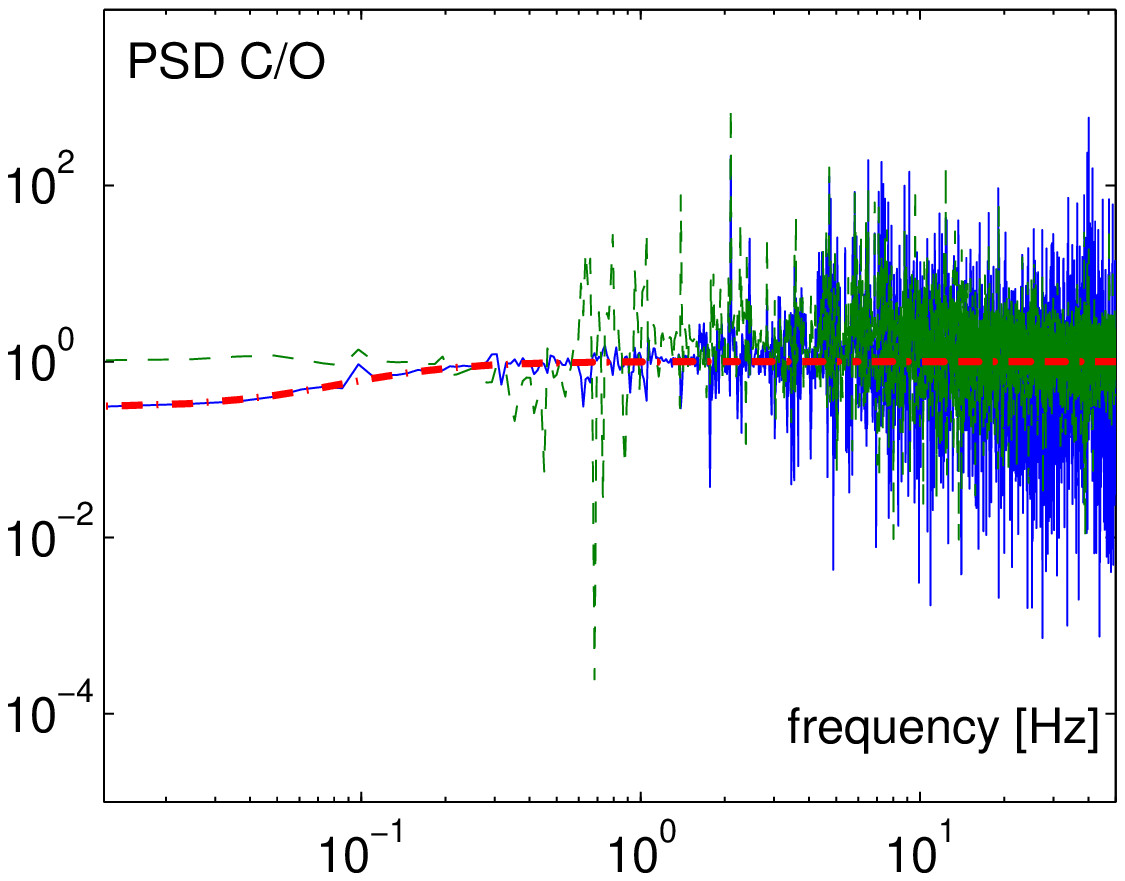}}	
	\put(0,40){\textbf{A.}}
	\put(65,40){\textbf{B.}}
	\put(130,40){\textbf{C.}}
	\put(136,40){Wavefront error after AO}
	\put(132,33){centralized}
	\put(132,26){distributed}
	\put(132,19){uniform gain}
	\put(132,12){open-loop}
	\put(159,33){2.8 nm}
	\put(159,26){2.8 nm}
	\put(159,19){20 nm}
	\put(159,12){16 nm}
	\put(155,8.5){{\scriptsize AO inoperative}}
	\put(152.5,10){{\Large (}}
	\put(179.5,10){{\Large )}}	
\end{picture}
\caption{\textbf{A}. Ratios of power spectral densities in the wavefront error associated to deformation mode $i=1000$, for a comprehensive EELT model subject to representative wind and sensor noise: blue solid curve [centralized control]/[open loop]; green dashed curve [distributed control]/[open loop]. Red dashdot curve represents $\vert (S(s))_{i,i} \vert^2 = \vert \tfrac{s+A_I(i)}{s+A_I(i)+K_I(i) \sqrt{\lambda_i}} \vert^2$, with $A_I(i)=0$ and $K_I(i)\sqrt{\lambda_i} = \tfrac{14.4}{2\pi}$~Hz from modal tuning. \textbf{B}. Same plot for deformation mode $j=2900$, $A_I(j)=\tfrac{0.39}{2\pi}$~Hz and $K_I(j)\sqrt{\lambda_j} = \tfrac{0.32}{2\pi}$~Hz. \textbf{C}. Overall root-mean-square wavefront error in the simulations, including a filter that models adaptive optics (AO) corrections. This effective filter is fictitious in open-loop, as high spatial frequency errors make AO unusable.}\label{fig:simu}
\end{figure*}

Proper control design does improve however the performance of the feedback system --- {\asedit within the limits imposed by distributed sensing}. The stability of the full simulation model, under a controller designed from simplified models (modal, LTSI), indicates the presence of sufficient margins for dynamic model uncertainties \textbf{[Prop.3]}.
We can estimate simulated disturbance rejection by taking the ratio of two simulation results, one with controller and one in open loop, with the same disturbance (wind) input. Figure \ref{fig:simu} shows the result for two representative deformation modes. The curves are in good agreement with the predicted $\vert (S(s))_{i,i} \vert^2$ under centralized control. (The irregularity of the graph at high frequencies just reflects that the denominator (disturbance input) of the ratio is nearly zero.) For the large-scale deformation mode 2900, the distributed controller performs slightly worse than the centralized one; the latter already does not reject much, due to its large $A_I(i)$ and low $K_I(i) \sqrt{\lambda_i}$. Overall, both controllers achieve similarly good rejection on small spatial scales and bad rejection on large spatial scales. In the EELT implementation, adaptive optics is added to the distributed sensing-based controller in order to reject the strong low spatial frequency disturbances and reach the overall $10$~nm precision requirement.

\begin{figure}[!t] 
    \centering
    \includegraphics[bb=20 180 585 610,clip=true,width=50mm]{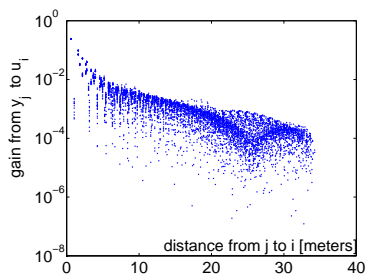}
    \caption{Gain given to sensors $j=1...5604$ in the centralized control input of segment $i$ (fixed to a typical value), as a function of distance to segment $i$.}
    \label{fig:locality}
\end{figure}

The simulations confirm that, under distributed sensing, the larger freedom offered by centralized control --- information transmission over large distances --- does not allow for significant improvement over the distributed controller. In agreement with similar theoretical investigation in \cite{Dahleh-Bamieh}, properly designed centralized controller gains are naturally ``distributed'': the gain from measurement at sensor $j$ to input at segment $i$ quickly decreases with the distance between $j$ and $i$, see Fig.~\ref{fig:locality}. 

\noindent\textbf{Remarks:} Spatiotemporal correlations also show up in model approximations. To avoid exciting actuator-backstructure dynamics at $x\un{Hz}$, a rule of thumb is to impose a closed-loop bandwidth below $\tfrac{x}{10}\un{Hz}$. From a finite-element model of the full telescope, low temporal frequency vibration modes also have low \emph{spatial} frequency. As a result, vibration modes significantly coupled to deformation modes $i=1$ to 2692 have frequencies above $\sim 30\un{Hz}$, allowing $3\un{Hz}$ control bandwidth. Only deformation modes $i>2693$ are coupled to lower temporal frequency vibrations  --- down to $\sim 3\un{Hz}$ --- and should thus be restricted to $0.3\un{Hz}$ bandwidth. It is remarkable that these bandwidth restrictions coincide with restrictions imposed for noise limitation.


\section{Conclusion}\label{sec:conclusion} 

This paper highlights how important limitations for linear control of distributed plants arise directly from \emph{local relative sensing} technology, even in the absence of any (e.g.~communication, decentralization) restrictions on the controller. We specifically show how distributed sensing can severely degrade the tradeoff between disturbance rejection and robustness to noise and to model uncertainties. The viewpoint allows correlations among subsystems and among disturbances. We illustrate our developments on the European Extremely Large Telescope's segmented primary mirror controller. The simulations suggest that with distributed sensing, a ``distributed controller'', computing each action from local measurements, achieves performance sensibly equal to an unrestricted centralized controller. Given distributed controllers' apparent adequacy for locally coupled distributed systems (see also \cite{Dahleh-Bamieh}), an exact study of the interplay between distributed sensing and structural limitations on the controller (distributed, leader-follower,...) could be interesting for future work. The role of \emph{relative} sensing in this context has recently attracted attention, see e.g.~\cite{Sarlette2009a,B2Ptac11}. ``Boundary conditions'' reminiscent of PDEs can also play an important role~\cite{Meurer11,Langbort2005}.

The paper more generally highlights the role of different sensor types for controlling distributed systems. A controller that relies on absolute measurements (``centralized sensing'') needs a common physical reference and achieves control performance uniform over all spatial frequencies. In multi-scale problems, the measurement accuracy attainable in this way may be insufficient and localized relative sensing is then used. But this seems to follow a waterbed effect: it improves control on spatial small-range signals at the expense of degradation on long-range signals. Centralized and distributed sensing are thus different spatial filters whose adequacy depend on the objective. The final EELT performance is obtained by combining the control loop of Section \ref{sec:EELT}, based on local edge sensors, and an adaptive optics controller based on a centralized wavefront-sensing. This suggests to combine sensors on a hierarchy of spatial scales for challenging distributed systems applications. 


\section*{Acknowledgments}

The authors thank C.~Bastin at U.Li\`ege and M.~Dimmler, B.~Sedghi, T.~Erm, B.~Bauvir from ESO for joint work on the EELT application. P.~Kokotovic is acknowledged for insightful comments on an early version of the manuscript. This paper presents research results of the Belgian Network DYSCO (Dynamical Systems, Control, and Optimization), funded by the Interuniversity Attraction Poles Program, initiated by the Belgian State, Science Policy Office. AS was an FNRS postdoctoral researcher at U.Li\`ege during part of this work.

	
\bibliographystyle{plain}
\bibliography{eso2.bib}


\section*{Appendix: design of the controllers}

\paragraph*{Centralized controller} Modulo approximating the high-frequency actuator resonances, both $G(s)$ and $B$ are diagonal (concatenated with appropriate all-zero matrices) in the modal basis $\modal{y} = Q^T y,\,Q^T u,\, U^T z$, so taking $K_I$ and $A_I$ diagonal in this same basis seems justified. The unobservable modes cannot be controlled, so $K_I(k) = 0$ for $k=2949$ to $2952$. The remaining components $1,2,...,N_0=2948$ are tuned as follows.
\vspace{2mm}\newline
First we discard leakage. For low $k$, the noise is negligible in Eq.~\eqref{eq:Modes}. Disturbances are divided by a given value on all modes by taking the same $(K)_{k,k} = K_0$, requiring
\begin{equation}\label{eq:Kival}
K_I(k) = K_0 / \sqrt{\lambda_k} \,, k = 1,2,...,2692\,.
\end{equation}
That is, less observable modes need larger $K_I(k)$. For large $k$, noise amplification makes \eqref{eq:Kival} a poor choice. We therefore must impose $(K)_{k,k} = K_1 \sqrt{\lambda_k}$ for $k>2692$ (where $\sqrt{\lambda_k} < 10 \sqrt{\lambda_1}$), i.e.~decreasing closed-loop gain, which translates into $K_I(k) = K_1$ constant. 
Our simulations take $K_0= 14.4$~rad/s, $K_1=5.7$~rad/s.
\vspace{2mm}\newline
Then we tune $A_I(k)$ with the approximate procedure around Proposition 2: for each $k$, we build a `bad' $\Delta$ according to \eqref{eq:ZeMethod1}, and assume $\Phi \approx \Phi_k$ as in \eqref{eq:ApproxPhi}. Writing the low-frequency approximation $y = G\cdot K\cdot z \approx 1 \cdot \tfrac{K_I}{s+A_I} \cdot z$ in time domain and replacing $z$ yields
$\tfrac{d}{dt} (\modal{y})_k = - \left(K_I(k)\sqrt{\lambda_k}(1+\phi_k)+A_I(k) \right)\, (\modal{y})_k \, =: - \mu_k\; (\modal{y})_k \,.$
We hence impose $A_I(k) = \max\{p_0 - K_I(k)\sqrt{\lambda_k}(1+\phi_k)\}$ to get $-\mu_k \leq -p_0 = -0.1$~Hz. 
This results in nonzero leakage for $k\geq 2846$, with maximum $A_I(N_0) = 0.7$~rad/s.
A numerical computation confirms that the exact closed-loop eigenvalues with this $\Delta$ are below $-p_0$.\vspace{-5mm}

\paragraph*{Distributed controller} To gain more insight, we factorize the controller into two steps. First each segment computes a local estimate $\hat{y}_k  \in \mathbb{R}^3$ of its own configuration from the ES measurements at its boundaries. Then a local controller computes $(u)_k  \in \mathbb{R}^3$ from the estimates $\hat{y}_j$ of $k$ and of its immediate neighbors. Thus $(u)_k$ depends only on measurements made at the edges of segment $k$ and of its immediate neighbors.
\vspace{2mm}\newline
Controller structure, plant, and measurement operator are all invariant from one segment to another, except for the mirror's boundary. Our spatial invariance assumption (i) approximates the mirror as infinite or periodic and (ii) imposes spatially invariant controller parameters. The system then decouples into one $3$-dimensional (piston,tip,tilt)-system per $2$-dimensional spatial frequency. We adapt the LTSI method of \cite{Dahleh-Bamieh,Stein2005,Gorinevsky2008} from SISO- to MIMO-subsystems in a straightforward way.
\vspace{2mm}\newline
The local estimation $\hat{y}_k$ is computed by pseudo-inverting, at each segment, the $12\times 3$ local matrix linking the 12 sensors at its edges to its own 3 degrees of freedom only. This amounts to segment $k$ assuming that its neighbors are perfectly positioned and only itself is displaced. That is of course not correct, so the local estimator acts like a first spatial filter\footnote{Other spatial filters can be envisioned --- central estimation is one, for other options see \cite{MacMyn2005} --- but none will alleviate the limitations imposed by distributed sensing.}, reflecting the difficulties of distributed sensing described with Fig.~\ref{fig:locestmodes}.
As a consequence, the chosen estimation filter will\vspace{-4mm}
\begin{itemize}\tightlist
\item[-] couple piston, tip and tilt (i.e.~it estimates a tilt component in $\hat{y}$ for a pure piston deformation $y$);
\item[-] underestimate piston (resp.~tip,tilt) dominated deformations of low (resp.~high) spatial frequency;
\item[-] slightly overestimate piston (tip,tilt) dominated deformations of high (low) spatial frequency.
\end{itemize}\vspace{-4mm}
We denote the output-to-estimation transfer function $(\hat{y})_{\xi} = H^*(\xi,s)\, (y)_{\xi}$, such that $H^*(\xi,s)$ is a $3 \times 3$ transfer function matrix for each $\xi$ and $s$. Our static plant approximation makes $H^*$ static, so it is trivially spatiotemporally factorized, facilitating tuning (see \cite{Gorinevsky2008}).
\vspace{2mm}\newline
Restricting the controller \eqref{eq:GeneralController} to local coupling, spatial invariance, and symmetry w.r.t.~axes reversal, we get\vspace{2mm}
\newline $
K_I^*(\xi) = k_{\alpha} + k_{\beta} \cos(\xi_1) + k_{\gamma} \cos(\xi_2) + k_{\delta} \cos(\xi_1\text{-}\xi_2)$
\newline $	
A_I^*(\xi) = a_{\alpha} + a_{\beta} \cos(\xi_1) + a_{\gamma} \cos(\xi_2) + a_{\delta} \cos(\xi_1\text{-}\xi_2)$\vspace{2mm}
\newline
where $k_i, \, a_i$, $i \in \{\alpha,\beta,\gamma,\delta \}$, are $3\times3$ matrices to tune. The LTSI tuning is performed in two steps. First, we approximate that $H^*$ is diagonal, i.e.~it leaves piston, tip and tilt uncoupled. This allows to directly apply the SISO tuning method of \cite{Stein2005}, using a large but linear optimization setting that minimizes leakage under robust stability and disturbance rejection constraints. We set the latter according to expected $d$ and $n$ on spatiotemporal frequencies and impose robustness as for centralized control: closed-loop eigenvalues $\leq -p_0 = -0.1$~Hz with any model uncertainties $\Delta$. The spatial invariance approximation requires to include extra robustness, so we increase $\varepsilon$ to $2.5\%$. In a second step, we add a stability constraint for the nominal MIMO system ($H^*$ not diagonal, $\Delta=0$). This yields a set of nonlinear constraints. The resulting nonlinear optimization is solved locally, starting from the solution of the first step. 
Simulating the resulting controller on the exact system confirms that intended properties hold despite the LTSI approximation in its design. Thus the system with less than $40$ segments diameter and a hole in the middle (see Fig.~\ref{fig:M1}) can be approximated as spatially invariant for controller design when its coupling is sufficiently local.		
\end{document}